\documentclass[12pt]{article}
\usepackage{amsmath}
\usepackage{amssymb}
\usepackage{amscd}

\renewcommand{\P}{\mathbb{P}}
\newcommand{\C}{\mathbb{C}}

\newcommand{\R}{\mathbb{R}}

\newcommand{\Z}{\mathbb{Z}}

\newcommand{\h}{{\rm h}}

\renewcommand{\H}{{\cal H}}

\newcommand{\cad}{{\it c.-\`a-d. }}

\newcommand{\lov}{{\rm lov}}
\newcommand{\vol}{{\rm vol}}

\newcommand{\Lone}{{{\rm L}^1}}

\newcommand{\dist}{{\rm dist}}

\newcommand{\FS}{{\rm FS}}

\newcommand{\voir}{{\it voir }}

\newenvironment{preuve}{\begin{trivlist} \item[]{\it Preuve---}}
{\par\hfill $\square$\end{trivlist}}
\title{Une borne sup\'erieure pour
l'entropie topologique  d'une application rationnelle}
\author{Tien-Cuong Dinh et Nessim Sibony}
\date{}
\begin{document}
\maketitle
\begin{abstract} Let $X$ be a complex projective manifold and $f$
a dominating rational map from $X$ onto $X$. We show that the
topological entropy $\h(f)$ of $f$ is bounded from above by the
logarithm of its
maximal dynamical degree.
\end{abstract}
Soit $X$ une vari\'et\'e projective complexe
de dimension $k\geq
2$, munie d'une forme de K\"ahler $\omega$ normalis\'ee par
$\int_X\omega^k=1$.
Soit
$f:X\longrightarrow X$ une application rationnelle dominante,
\cad localement ouverte en un point g\'en\'erique de
$X$. Notons $\lambda_l(f)$ le degr\'e dynamique d'ordre $l$ de
$f$, $1\leq l\leq k$, et $\h(f)$ l'entropie
topologique de $f$. Ils sont d\'efinis plus loin.
Il s'agit de montrer que $\h(f)\leq \max_{1\leq l\leq k}
\log\lambda_l(f)$. Un interm\'ediaire utile est
$\lov(f)$, c'est un indicateur 
de la croissance du volume des graphes des it\'er\'es
de $f$. Plus pr\'ecis\'ement, soit $\Gamma_n$ le graphe de
l'application $(f,f^2,\ldots,f^n)$ o\`u 
$f^i:=f\circ\cdots\circ f$ ($i$ fois). On pose
$$\lov(f):=\limsup_{n\rightarrow \infty}
\frac{1}{n}\log\vol(\Gamma_n).$$
\par
Utilisant une
in\'egalit\'e de Lelong, Gromov \cite{Gromov1}
a montr\'e que
$\h(f)\leq \lov(f)$ lorsque $f$ est
holomorphe. Sa preuve reste
valable pour les applications rationnelles.
Nous allons montrer que
$\lov(f) = \max_{1\leq l\leq k}\log\lambda_l(f)$ et en d\'eduire
le th\'eor\`eme suivant.
\par
\ 
\\
{\bf Th\'eor\`eme 1} {\it
Soit $X$ une vari\'et\'e projective complexe 
de dimension
$k\geq 2$ et soit $f:X\longrightarrow X$ une application rationnelle
dominante. Alors 
$$\h(f)\leq \lov(f)=\max_{1\leq l\leq k} \log\lambda_l(f).$$}
\par
\
\par
Ce r\'esultat a \'et\'e annonc\'e dans \cite{Friedland}.
L'auteur a utilis\'e  l'in\'egalit\'e $(f\circ
f)^*\leq f^*\circ f^*$ \cite[lemma 3]{Friedland}
qui n'est pas valable lorsque ces
op\'erateurs agissent sur les cycles analytiques. Un
contre-exemple est donn\'e dans \cite{Guedj}. 
Dans notre approche, on ne d\'efinit pas $f^*$ sur tous les
courants ni m\^eme sur les cycles. Le produit des courants n'est
consid\'er\'e que l\`a o\`u ils sont lisses. 
\par
Rappelons quelques notions 
(\voir par exemple \cite{Bowen,Sibony2,Guedj, DinhSibony1}).
Notons $\Gamma$ le graphe de $f$ dans $X\times X$. C'est un
sous-ensemble analytique irr\'eductible de dimension $k$.
Lelong a montr\'e que l'int\'egration
sur la
partie r\'eguli\`ere de $\Gamma$ d\'efinit un courant positif
ferm\'e $[\Gamma]$ de bidimension $(k,k)$.
Soient $\pi_1$ et
$\pi_2$ les projections de $X\times X$ sur le premier et le
second facteur. Pour tout ensemble $Y\subset X$, posons
$$f(Y):=\pi_2(\pi_1^{-1}(Y)\cap \Gamma) \ \mbox{ et }\  
f^{-1}(Y):=\pi_1(\pi_2^{-1}(Y)\cap \Gamma).$$
\par
Notons $I_f$ l'ensemble des points $x\in X$ tels que
$\dim \pi_1^{-1}(x)\cap\Gamma\geq 1$. C'est {\it
l'ensemble des points
d'ind\'etermination de $f$}. Il est de codimension au moins 2.
On a $\dim \pi_1^{-1}(I_f)\cap\Gamma \leq k-1$.
Pour toute forme $\varphi$ lisse de bidegr\'e $(l,l)$ posons
\begin{eqnarray}
f^*(\varphi) & := & (\pi_1)_*(\pi_2^*(\varphi)\wedge [\Gamma])
\end{eqnarray}
Le courant $f^*(\varphi)$
est lisse sur $X\setminus I_f$. Puisqu'il est de masse finie et sans masse 
sur $I_f$, ses
coefficients sont
dans $\Lone$.
En particulier, $f^*(\varphi)$ ne
charge pas les sous-ensembles analytiques propres de $X$.
L'op\'erateur $f^*$ est continu de l'espace des formes lisses
dans l'espace des formes \`a coefficients dans $\Lone$.
\par
Posons $\delta_l(f):=\int_X f^*(\omega^l)\wedge
\omega^{k-l}$ pour $1\leq l\leq k$.
On d\'efinit {\it le degr\'e dynamique d'ordre $l$} de
$f$ par
$$\lambda_l(f):=\limsup_{n\rightarrow\infty} [\delta_l(f^n)]^{1/n}.$$
On verra que la suite $[\delta_l(f^n)]^{1/n}$
est en  toujours convergente (corollaire 7).
{\it Le degr\'e topologique} $d_t:=\lambda_k(f)$ est \'egal au nombre
de pr\'eimages par $f$ d'un point g\'en\'erique de $X$.
\par
Soit
$$\Omega_f:=X\setminus \cup_{n\in\Z} f^n(I_f).$$
C'est un ensemble invariant par $f$ et $f^{-1}$. On dira qu'une
famille $F\subset\Omega_f$ est {\it $(n,\epsilon)$-s\'epar\'ee},
$\epsilon>0$, si
$$\max_{0\leq i\leq n-1} \dist(f^i(x),f^i(y))\geq \epsilon
\ \mbox{ pour }\ x,y\in F \mbox{ distincts}.$$
{\it L'entropie topologique} (\voir \cite{Bowen})
$\h(f)$ est d\'efinie par
$$\h(f):=\sup_{\epsilon>0}\left(\limsup_{n\rightarrow\infty}
\frac{1}{n}
\log \max\big\{\#F, \ F\
(n,\epsilon)\mbox{-s\'epar\'ee}\big\}\right).$$
\par
Notons $\Gamma_n$ l'adh\'erence dans $X^n$ de l'ensemble des points
$$(x,f(x),\ldots, f^{n-1}(x)),\ \ \ x\in\Omega_f.$$ C'est un
sous-ensemble analytique de dimension $k$ de $X^n$.
Soient $\Pi_i$ les projections de $X^n$ sur ses facteurs. On
munit $X^n$ de la forme de K\"ahler $\omega_n:=\sum
\Pi_i^*(\omega)$. On a
$$\lov(f):=\limsup_{n\rightarrow\infty} \frac{1}{n}
\log(\vol(\Gamma_n)):=
\limsup_{n\rightarrow\infty}\frac{1}{n}\log
\left(\int_{\Gamma_n}\omega_n^k\right).$$
On verra plus loin que la suite 
$\big(\frac{1}{n}\log\int_{\Gamma_n}\omega_n^k\big)$ est toujours convergente.
Utilisant une in\'egalit\'e de Lelong,
Gromov \cite{Gromov1}
a montr\'e que $\h(f)\leq \lov(f)$. Dans la
suite, nous montrons que $\lov(f)=\max \log \lambda_l(f)$.
\par
Pour tout courant positif ferm\'e de bidegr\'e $(l,l)$ sur $X$,
notons $\|S\|:=\int_X S\wedge \omega^{k-l}$ {\it la masse} de $S$.
Puisqu'on a suppos\'e $\int_X\omega^k=1$, les courants $\omega^l$
sont de masse 1. Pour les r\'esultats fondamentaux sur les
courants positifs ferm\'es nous renvoyons \`a Lelong
\cite{Lelong} et Demailly \cite{Demailly}.
Notre outil principal est le lemme suivant. 
\\
\
\\
{\bf Lemme 2.} {\it
Il existe une constante $c>0$, qui ne d\'epend que de
$X$, telle que pour tout courant positif ferm\'e $S$ de bidegr\'e
$(l,l)$ sur $X$ on puisse trouver une suite de courants positifs
ferm\'es lisses $(S_m)_{m\geq 1}$, de bidegr\'e $(l,l)$,
v\'erifiant les propri\'et\'es suivantes
\begin{enumerate}
\item La suite $(S_m)$ converge vers un courant positif ferm\'e $S'$.
\item $S'\geq S$, c'est-\`a-dire que le courant $S'-S$ est positif.
\item On a $\|S_m\|\leq c\|S\|$ pour tout $m\geq 1$.
\end{enumerate}
}
\par
\begin{preuve}
Soit $\omega_\FS$ la forme de Fubini-Study de
$\P^k$ normalis\'ee par $\int_{\P^k}\omega_\FS^k=1$.
Rappelons que les groupes de cohomologie de Dolbeault
$\H^{l,l}(\P^k,\R)$ sont de dimension $1$. En particulier, tout
courant positif ferm\'e
$R$ de bidegr\'e $(l,l)$ sur $\P^k$ est cohomologue \`a
$\|R\|\omega_\FS^l$. 
\par
Puisque $X$ est projective,
on peut choisir une famille finie d'applications
holomorphes surjectives $\Psi_i$, $1\leq i\leq s$, 
de $X$ dans $\P^k$ telles qu'en tout
point $x\in X$ au
moins l'une des applications $\Psi_i$ soit de rang maximal.
Il suffit de plonger $X$ dans un $\P^N$ et de prendre une famille
de projections sur $\P^k$.
Posons $T_i:=(\Psi_i)_*(S)$.
L'op\'erateur $(\Psi_i)_*$ \'etant continu, il
existe une constante $c_1>0$ ind\'ependante de $S$
telle que $\|T_i\|\leq c_1\|S\|$.
\par
La vari\'et\'e $\P^k$ \'etant homog\`ene, 
si $R$ est un courant positif
ferm\'e dans $\P^k$, il existe des courants
positifs ferm\'es lisses $(R_m)$ tendant vers $R$ avec
$\|R_m\|=\|R\|$. 
Donc il existe 
des courants positifs ferm\'es lisses $T_{i,m}$ de
bidegr\'e $(l,l)$ sur $\P^k$ qui convergent
faiblement vers $T_i$ et qui v\'erifient  $\|T_{i,m}\|=\|T_i\|$.
On a donc $\|T_{i,m}\|\leq c_1\|S\|$.
\par
Posons $S_m:=\sum_{i=1}^s (\Psi_i)^*(T_{i,m})$.
Estimons 
la masse de
$(\Psi_i)^*(T_{i,m})$:
\begin{eqnarray*}
\|(\Psi_i)^*(T_{i,m})\| & = & \int_X (\Psi_i)^*(T_{i,m})\wedge
\omega^{k-l}\\
& = & \|T_{i,m}\|\int_X (\Psi_i)^*(\omega_\FS^l)
\wedge\omega^{k-l}\\
& \leq & c_1\|S\|\int_X (\Psi_i)^*(\omega_\FS^l)
\wedge\omega^{k-l}\\
& \leq & c_2\|S\|
\end{eqnarray*}
pour une constante $c_2>0$
ind\'ependante de $S$. 
Donc la masse de $S_m$ est major\'ee par $c\|S\|$ avec $c:=sc_2$.
Quitte \`a extraire une sous-suite, on peut supposer que la suite
$(S_m)$ tend faiblement vers un courant $S'$.
Au voisinage de chaque point $x\in X$, on v\'erifie,
puisque l'un des $\Psi_i$ est un biholomorphisme
local, que $S'-S$ est positif. On peut bien s\^ur
choisir les constantes $c_1$, $c_2$ et $c$ ind\'ependantes de
$l$, $1\leq l\leq k$.  
\end{preuve}
{\bf Remarque 3.} 
L'ensemble des classes de courants positif 
ferm\'es de bidegr\'e $(l,l)$ 
et de masse 1  est born\'e dans $\H^{l,l}(X,\R)$. 
Il existe  donc une
constante $\alpha_X>0$ telle que la classe de
$\alpha_X\omega^l -T$ soit
repr\'esent\'ee par une forme lisse positive pour tout courant
positif ferm\'e $T$ de bidegr\'e $(l,l)$ et de masse plus petite ou
\'egale \` a 1.
On dira que $T$ est {\it cohomologiquement domin\'e}
par $\alpha_X\omega^l$.
La propri\'et\'e ci-dessus est valable pour toute vari\'et\'e k\"ahl\'erienne
compacte. Dans le lemme 2,
les courants $S_m$  sont cohomologiquement domin\'es
par $c_X\|S\|\omega^l$ o\`u $c_X:=c\alpha_X$.
\\
\
\par
Notons ${\cal C}_f$ l'ensemble des points au voisinage desquels
$f$ n'est pas une application holomorphe localement inversible. 
Posons 
$\Omega_{1,f}:=
X\setminus {\cal C}_f$. C'est un
ouvert de Zariski de $X$.
Nous allons d\'efinir $f^*(S)$ lorsque $S$ est un courant positif ferm\'e
de bidegr\'e $(l,l)$ sur $X$.
Le courant $f^*(S)$ est bien d\'efini sur
$\Omega_{1,f}$. Si sa masse sur $\Omega_{1,f}$, qui est d\'efinie par
$\|f^*(S)\|:=\int_{\Omega_{1,f}} f^*(S)\wedge
\omega^{k-l}$, est finie, d'apr\`es Skoda \cite{Skoda},
son prolongement trivial $\widetilde{f^*(S)}$ est un courant
positif ferm\'e sur $X$. Le lemme suivant montre que c'est le cas. Nous utilisons 
par la suite cette extension par 0.
\\
\
\\
{\bf Lemme 4.} {\it
Soit $S$ un courant positif ferm\'e de bidegr\'e $(l,l)$ sur $X$.
Alors $\|f^*(S)\|\leq c_X\delta_l(f)\|S\|$. En particulier,
$\widetilde{f^*(S)}$ est positif ferm\'e dans $X$ et sa masse est
born\'ee par $c_X\delta_l(f)\|S\|$.
}
\begin{preuve}
Soit $(S_m)$ la suite de courants lisses
v\'erifiant le lemme 2 (appliqu\'e au courant $S$).
D'apr\`es la remarque 3, ces courants sont
cohomologiquement domin\'es par $c_X\|S\|\omega^l$.
On en d\'eduit que la masse de
$f^*(S_m)$, qui se calcule cohomologiquement,
est major\'ee par
$c_X \delta_l(f)\|S\|$. Plus pr\'ecis\'ement, on a pour tout compact 
$K\subset\Omega_{1,f}$  
\begin{eqnarray*}
\int_K f^*(S)\wedge \omega^{k-l} & \leq &
\int_K f^*(S')\wedge \omega^{k-l}\leq
\lim_{m\rightarrow\infty} \int_X f^*(S_m)\wedge
\omega^{k-l} \\
& \leq & c_X\|S\|\int_X
f^*(\omega^l)\wedge\omega^{k-l}
=  c_X\delta_l(f)\|S\|.
\end{eqnarray*}
Ceci implique le lemme.
\end{preuve}
\
\\
{\bf Lemme 5.} {\it Soit $\epsilon>0$. Il existe une constante
$c_\epsilon>0$ telle qu'on ait
$$\int_{\Omega_f} (f^{n_1})^*\omega\wedge \ldots\wedge
(f^{n_k})^*\omega \leq c_\epsilon
(\max_{1\leq l\leq k}\lambda_l+\epsilon)^{n_1}$$
pour tous les entiers naturels $n_1,\ldots,
n_k$ v\'erifiant $n_1\geq\cdots\geq n_k\geq 0$. }
\
\\
\begin{preuve} Posons $\lambda_\epsilon:=\max_{1\leq l\leq
k}\lambda_l+\epsilon$. Soit $c>0$ une constante
telle que $\delta_l(f^n)\leq c\lambda_\epsilon^n$ pour
tout $n\geq 0$ et pour tout $l$ avec $1\leq l\leq k$. Soit
$\Omega_{n,f}:=X\setminus \cup_{0\leq i\leq n-1} f^{-i}({\cal C}_f)$.
C'est un ouvert de Zariski de $X$. 
Montrons par r\'ecurrence sur $s$, $0\leq s\leq k$, que 
pour tous $n_1,\ldots,
n_s$ v\'erifiant $n_1\geq\cdots\geq n_s\geq 0$ on a
$\|T_s\|\leq c^s c_X^s \lambda_\epsilon^{n_1}$,
$c_X$ \'etant la constante de la remarque 3 et 
$$T_s:=(f^{n_1})^*\omega\wedge \ldots\wedge (f^{n_s})^*\omega,
\ \ \ \ T_0:=1.$$
\par
C'est clair au rang $s=0$. Supposons le au rang $s-1$,
$1\leq s\leq k$. On a
$\|T'_{s-1}\|\leq c^{s-1} c_X^{s-1} \lambda_\epsilon^{n_1-n_s}$
o\`u 
$$T'_{s-1}:=(f^{n_1-n_s})^*\omega\wedge
\ldots\wedge (f^{n_{s-1}-n_s})^*\omega.$$
Le courant $T'_{s-1}$ \'etant de masse finie sur
$\Omega_{n_1-n_s,f}$, d'apr\`es le th\'eor\`eme de Skoda
\cite{Skoda}, son prolongement trivial $\widetilde{T_{s-1}'}$
dans $X$ est un courant positif
ferm\'e dont la masse est major\'ee par
$c^{s-1} c_X^{s-1} \lambda_\epsilon^{n_1-n_s}$.
Utilisant le lemme 4 appliqu\'e au courant
$S=\widetilde{T'_{s-1}}\wedge\omega$
et \`a l'application $f^{n_s}$, on obtient
$$\|T_s\|=\|(f^{n_s})^*(T'_{s-1}\wedge\omega)\|\leq c_X
\delta_l(f^{n_s}) \|T'_{s-1}\|\leq
c^sc_X^s\lambda_\epsilon^{n_1}.$$  
Ceci termine la r\'ecurrence. Pour $s=k$, on
obtient le lemme avec $c_\epsilon:=c^kc_X^k$.
\end{preuve}
{\bf Fin de la d\'emonstration du th\'eor\`eme 1.} Il nous faut
estimer $\lov(f)$. D'apr\`es le lemme 5, on a pour tout
$\epsilon>0$ 
\begin{eqnarray*}
\vol(\Gamma_n) & = & \sum_{0\leq i_1,\ldots i_k\leq n-1}
\int_{\Omega_f}
(f^{i_1})^*\omega \wedge \ldots \wedge  (f^{i_k})^*\omega\\
& \leq & c_\epsilon n^k (\max_{1\leq l\leq
k}\lambda_l(f)+\epsilon)^{n-1}.
\end{eqnarray*}
D'o\`u $\lov(f)\leq \max_{1\leq l\leq k} \log\lambda_l(f)$.
On a aussi $\lov(f)\geq \max_{1\leq l\leq k} \log\lambda_l(f)$
car
$$\vol(\Gamma_n)\geq \int_{\Omega_f}
(f^{n-1})^*\omega^l \wedge \omega^{k-l} =\delta_l(f^{n-1}).$$
\par
\hfill $\square$
\\
\
\\
{\bf Proposition 6.} {\it Soient $f$ et $g$ deux applications
rationnelles de $X$ dans $X$. On a
$$\delta_l(f\circ g)\leq c_X\delta_l(f)\delta_l(g).$$
}
\begin{preuve}
Par d\'efinition de $\delta_l(f)$, on a $\|f^*\omega^l\|=\delta_l(f)$.
Le courant $(f\circ g)^*\omega^l$, qui ne charge pas les
sous-ensembles analytiques propres de $X$, est \'egal \`a
$g^*(\widetilde{f^*\omega^l})$ sur
$\Omega_{1,g}\cap\Omega_{1,f\circ g}$. 
Le lemme 4, appliqu\'e \`a $g$ et au
courant $\widetilde{f^*\omega^l}$, entra\^{\i}ne que
$$\delta_l(f\circ g)=\|g^*(\widetilde{f^*\omega^l})
\|\leq c_X\delta_l(g)\|f^*\omega^l\| =
c_X\delta_l(f)\delta_l(g).$$ 
\end{preuve}
\
\\
{\bf Corollaire 7.} {\it La suite $[\delta_l(f^n)]^{1/n}$
est convergente. Les degr\'es dynamiques $\lambda_l$ de $f$ sont
des invariants birationnels. 
}
\begin{preuve} D'apr\`es la proposition 6, on a
$\delta_l(f^{m+n}) \leq c_X \delta_l(f^m)\delta_l(f^n)$ pour
tous $m,n\geq 1$. Ceci implique que la suite $[\delta_l(f^n)]^{1/n}$
converge vers $\inf_{n\geq 1}[\delta_l(f^n)]^{1/n}$.
\par
Soit $g$ une application birationnelle de $X$ dans $X$. Posons
$h:=g\circ f\circ g^{-1}$.
On a 
$$\delta_l(h^n)=\delta_l(g\circ f^n\circ g^{-1})\leq c_X^2
\delta_l(g)\delta_l(g^{-1})\delta_l(f^n).$$
Donc $\lambda_l(f)\leq \lambda_l(h)$.
Puisque $f=g^{-1}\circ h\circ g$, on a aussi
$\lambda_l(h)\leq \lambda_l(f)$. 
\end{preuve}
\
\\
{\bf Remarques 8.} {\bf a.} Russakovskii et Shiffman
\cite{RussakovskiiShiffman} ont montr\'e l'in\'egalit\'e
$\delta_l(f\circ g)\leq \delta_l(f)\delta_l(g)$ lorsque
$X=\P^k$. Diller et Favre \cite{DillerFavre}
ont d\'ecrit pr\'ecis\'ement la
croissance de $\delta_1(f^l)$ lorsque $f$ est une application
bim\'eromorphe sur une surface complexe.
Dans le cas de dimension $k\leq 3$
et dans le cas des vari\'et\'es homog\`enes,
les r\'esultats ci-dessus ont \'et\'e d\'emontr\'es par
Vincent Guedj \cite{Guedj}.
Il a alors prouv\'e l'existence d'une unique mesure
invariante d'entropie maximale $\log d_t(f)$ pour toute
application rationnelle $f$ v\'erifiant $d_t(f)>
\lambda_l(f)$, $1\leq l\leq k-1$ (pour la m\'ethode \voir \'egalement 
Briend-Duval
\cite{BriendDuval} ainsi que \cite{DinhSibony1,DinhSibony3}). 
Le th\'eor\`eme 1 permet
d'\'etendre ce r\'esultat au cas d'une 
vari\'et\'e projective quelconque. 
\par
{\bf b.} D'apr\`es Iskovkikh-Manin \cite{IskovkikhManin}, il existe 
des vari\'et\'es $X$ lisses non rationnelles de dimension 3 
dans $\P^4$ qui sont unirationelles, \cad pour lesquelles
il existe une application rationnelle de 
rang maximal $f:\P^3\longrightarrow X$. On peut donc composer
une projection holomorphe $g$ de $X$ sur $\P^3$ avec $f$ pour obtenir 
beaucoup d'applications rationnelles d'entropie positive sur $X$.
Nous remercions F. Campana et N. Mok qui nous ont indiqu\'e ces exemples.
\par
Il est facile cependant de construire
des correspondances sur les vari\'et\'es projectives (\voir par
exemple \cite{ClozelUllmo,Voisin,Dinh,DinhSibony2}).
La proposition 6 et le
corollaire 7   
restent valables pour les correspondances.
Rappelons qu'une correspondance sur $X$
est la donn\'ee d'un ensemble
analytique $\Gamma\subset X\times X$ de dimension $k$
dont les images de chaque
composante par $\pi_1$ et $\pi_2$ sont \'egales \`a $X$. On peut
poser $f:=\pi_2\circ(\pi_{1|\Gamma})^{-1}$. L'image
r\'eciproque d'une forme lisse
est d\'efinie par l'\'equation (1). Les
degr\'es dynamiques sont d\'efinis de fa\c con analogue que pour les 
applications rationnelles. 
\par
{\bf c.} Soit $\pi:Y\longrightarrow X$ une application
holomorphe surjective o\`u $X$ et $Y$ sont des vari\'et\'es
projectives complexes.
Soit $T$ un courant positif ferm\'e sur $X$.
L'images r\'eciproque $\pi^*(T)$
de $T$ par $\pi$ est bien d\'efinie
sur un ouvert de Zariski de $Y$ (l\`a o\`u $f$ est une submersion
locale). Le lemme 2 permet de prolonger $\pi^*(T)$ en
courant positif ferm\'e $\widetilde{\pi^*(T)}$ dans $Y$.
L'op\'erateur $T\mapsto \widetilde{\pi^*(T)}$ est semi-continu
inf\'erieurement.
Plus pr\'ecis\'ement, si $T_n\rightarrow T$, $\lim
\widetilde{\pi^*(T_n)}\geq \widetilde{\pi^*T}$. 
Cette d\'efinition est utile dans
le cadre des courants dynamiques. Nous reviendrons sur cette
question dans un prochain travail.
Notons que M\'eo \cite{Meo}
a donn\'e un exemple qui montre qu'on ne peut pas toujours
d\'efinir $\pi^*(T)$ dans le cas o\`u $X$ et $Y$ ne
sont pas compactes. Il a aussi donn\'e une
d\'efinition de $\pi^*(T)$ dans le cas local et lorsque $\pi$ est une
application \`a fibres discr\`etes.
On ne sait pas si sa
d\'efinition est ind\'ependante des coordonn\'ees.
\small

Tien-Cuong Dinh et Nessim Sibony,\\
Math\'ematique - B\^at. 425, UMR 8628, 
Universit\'e Paris-Sud, 91405 Orsay, France. \\
E-mails: Tiencuong.Dinh@math.u-psud.fr et
Nessim.Sibony@math.u-psud.fr.
\end{document}